\documentclass[12pt]{amsart}

\newcommand{\haus}{\operatorname{\mathcal H}}

\newcommand{\Vol}{\operatorname{vol}}
\newcommand{\dist}{\operatorname{dist}}
\newcommand{\C}{{\mathbf C}}
\newcommand{\dbar}{\overline \partial}

\newtheorem*{theorem}{Theorem}
\newtheorem*{stability}{Stability lemma for the
Bergman kernel function}
\newtheorem*{barbell}{Barbell lemma}

\begin{document}

\title{The Lu Qi-Keng Conjecture Fails Generically}

\author{Harold P. Boas}

\address{Department of Mathematics, Texas A\&M
University, 
College Station,
TX 77843--3368}

\thanks{This research was partially
supported by NSF grant number DMS-9203514.\\
This is a revised version of a paper 
presented at meeting 895 of the
American Mathematical Society.}

\email{boas@@math.tamu.edu}


\maketitle

\begin{abstract}
The bounded domains of holomorphy in~$\C^n$ whose Bergman
kernel functions are zero-free form a nowhere
dense subset (with respect to a variant of the
Hausdorff distance) of all bounded domains of
holomorphy.
\end{abstract}

A domain in~$\C^n$ is called a {\it Lu
Qi-Keng domain\/} if its Bergman kernel function
has no zeroes. Lu Qi-Keng \cite{LuQiKeng} raised the question of
which domains, besides  the ball
and the polydisc, have this property. A motivation
for the question is that vanishing of the Bergman
kernel function obstructs the global definition of
Bergman representative coordinates. Over the years
since Lu Qi-Keng's paper appeared, 
various versions of a {\it Lu Qi-Keng
conjecture\/} have been mooted to the effect that all domains, or
most domains, or all domains satisfying some
geometrical hypothesis, are Lu Qi-Keng domains. 

In the complex plane~$\C^1$, a bounded domain with smooth
boundary is a Lu Qi-Keng domain if and only if it
is simply connected \cite{Suita Yamada} (and thus biholomorphically
equivalent to the disc). I have given a
counterexample \cite{Boas}
showing that no analogous topological
characterization of Lu Qi-Keng domains can hold in
higher dimensions: there exists (in~$\C^2$,
and similarly in~$\C^n$ for $n>2$) a bounded,
strongly pseudoconvex, contractible domain with
$C^\infty$~regular boundary whose Bergman kernel
function does have zeroes. 

In this note, I show that the Lu Qi-Keng domains
of holomorphy may be viewed as exceptional: they form
a nowhere dense set with respect to a suitable
topology. Thus, contrary to former expectations,
it is the normal situation for the Bergman kernel function of a
domain to have zeroes.

To formulate the result precisely, I need a metric
on bounded open sets. Since I impose no restriction on
the regularity of the boundaries of the sets, some
variant of the Hausdorff metric will be
appropriate. The Hausdorff distance~$\haus$ is
normally defined for nonempty, bounded, closed
sets by the property that 
$\haus(A,B)<\epsilon$ if and only if each point of~$A$
has Euclidean distance less than~$\epsilon$ from
some point of~$B$, and vice versa.

After the seminal paper of Ramadanov \cite{Ramadanov}, it is clear
in the context of the Bergman
kernel function that if a sequence of open sets~$\{\Omega_j\}$
is going to be said to converge to an open set~$\Omega$, then every
compact subset of~$\Omega$ should eventually be
contained in~$\Omega_j$. It is less clear what
requirement should be imposed if the
$\Omega_j$~contain points outside of~$\Omega$. 
The example \cite[p.~39]{Skwarczynski} 
\cite[p.~280]{Ramadanov 1983}
of decreasing concentric disks in the
complex plane converging to a disk with a slit
removed shows that it is inadequate to require
merely that for every
open neighborhood of the closure~$\overline\Omega$,
eventually $\Omega_j$ is contained in the
neighborhood. 

I shall consider two different notions of
convergence of open sets in~$\C^n$. 
Both have the property that if
$\Omega_j\to\Omega$, then the $\Omega_j$
eventually swallow every compact subset of~$\Omega$.
However, they differ in what they require about
the sets $\Omega_j\setminus\Omega$.

First I define a metric~$\rho_1$ on bounded,
nonempty, open sets via
$\rho_1(U,V)=\haus(\overline U,\overline V) +
\haus(\partial U,\partial V)$. If
$\rho_1(\Omega_j, \Omega)\to 0$, then
the~$\Omega_j$ eventually swallow every compact
subset of~$\Omega$ and are eventually swallowed by
every open neighborhood of~$\overline\Omega$. The
converse holds when $\Omega$ equals the interior of
its closure, but not in general.  By requiring
that both the closures and the boundaries
converge, convergence in the metric~$\rho_1$
eliminates examples like the one above involving
slits or punctures in the limit domain.

The metric~$\rho_1$ can also be thought of in
terms of functions. Define the distance
function~$d_U$ of an open set~$U$ via
$d_U(z)=\dist(z, \C^n\setminus U)$, where
$\dist$~denotes the Euclidean distance. Then 
$\Omega_j\to\Omega$  according to the
metric~$\rho_1$ if and only if the continuous
functions~$d_{\Omega_j}$ converge uniformly
on~$\C^n$ to~$d_\Omega$ 
and the $d_{\C^n\setminus\overline{\Omega_j}}$
converge uniformly to~$d_{\C^n\setminus\overline\Omega}$.

In some contexts---one will appear below---it is
useful to relax the hypothesis on how the
sets $\Omega_j\setminus\Omega$ behave. They
could be required to shrink in volume (Lebesgue
$2n$-dimensional measure), but not necessarily in
terms of Euclidean distance from~$\Omega$. I therefore
introduce a second metric~$\rho_2$ on bounded open
sets via $\rho_2(U,V)=\Vol(U\setminus
V)+\Vol(V\setminus
U)+\sup_{z\in\C^n}|d_U(z)-d_V(z)|$. 
Convergence of $\Omega_j$ to~$\Omega$ in this metric
allows $\Omega_j$ to have a long thin tail
whose width shrinks to zero but whose length does
not shrink. 

So far I have not assumed that the open sets in
question are connected. It is easy to see that the
Bergman kernel function $K(w,z)$ for a disconnected open
set is identically equal to zero if $w$ and~$z$
are in different connected components, while if
the points are in the same connected component,
then $K(w,z)$ equals the Bergman kernel function
of that component. I will say that a (possibly)
disconnected, bounded, nonempty, open set is a {\it Lu Qi-Keng open
set\/} if its Bergman kernel function has no
zeroes when the two variables are in the same
connected component.

\begin{theorem}
The Lu Qi-Keng open sets are nowhere dense in each of
the following metric spaces, where the metric is~$\rho_1$:
\begin{enumerate}
\item the bounded pseudoconvex open sets;
\item the bounded connected pseudoconvex open sets
(domains of holomorphy);
\item the bounded strongly pseudoconvex open sets;
\item the bounded connected strongly pseudoconvex
open sets.
\end{enumerate} 
If one considers only
open sets of Euclidean diameter less than some
fixed constant~$M$, then the same assertion holds
when the metric is taken to be~$\rho_2$.
\end{theorem}

I will base the proof of the theorem on the
following two folklore lemmas. 
The ideas of the proofs are all in the literature, but
since I do not know
a reference for precisely these formulations, I
will indicate proofs of the lemmas after the proof of the
theorem.

\begin{stability}
Let $\{\Omega_j\}$ be a sequence of 
bounded pseudoconvex open sets that converges,
in the sense of either $\rho_1$
or~$\rho_2$, to a nonempty bounded open
set~$\Omega$; in the case of~$\rho_2$, assume also
that the $\Omega_j$ have uniformly bounded
diameters (this is automatic in the case
of~$\rho_1$). 
Suppose $U$~is a connected component
of~$\Omega$ that has $C^\infty$~regular boundary
and that is separated from the rest of~$\Omega$
(that is, the closure of~$U$ is disjoint from the
closure of $\Omega\setminus U$). Then the Bergman
kernel functions of the~$\Omega_j$ converge to the
Bergman kernel function of~$U$ uniformly on
compact subsets of~$U\times U$.
\end{stability}

In the statement of the stability lemma,
pseudoconvexity of the limit set~$\Omega$ is
automatic: the limit function $-\log d_\Omega$ inherits
plurisubharmonicity from the 
functions $-\log d_{\Omega_j}$, which converge
uniformly on compact subsets of~$\Omega$.

A special case of considerable interest is when
the $\Omega_j$ and $\Omega$ are all bounded,
connected, pseudoconvex domains with $C^\infty$~regular
boundaries. The lemma then says that 
if the~$\Omega_j$ converge to~$\Omega$, in the
sense that $\Omega_j$ eventually swallows every
compact subset of~$\Omega$ and the volume of
$\Omega_j\setminus\Omega$ tends to zero, then
the Bergman kernel
functions of the~$\Omega_j$ converge uniformly on
compact subsets of $\Omega\times\Omega$ to the
Bergman kernel function of~$\Omega$.

The $C^\infty$~regularity hypothesis in the lemma
can be reduced to $C^2$~regularity, but I shall
not prove this here.

I take the name of the second lemma from 
\cite[Chap.~5, Exercise~21]{Krantz}.

\begin{barbell}
Suppose $G_1$ and~$G_2$ are bounded, connected,
strongly pseudoconvex domains in~$\C^n$ with
$C^\infty$~regular boundaries and with disjoint
closures. Let $\gamma$~be a smooth curve (that is,
a $C^\infty$~embedding of $[0,1]$ into~$\C^n$)
joining a boundary point of~$G_1$ to a boundary
point of~$G_2$, and otherwise disjoint from the
closures of $G_1$ and~$G_2$, and let $V$~be an arbitrary
neighborhood in~$\C^n$ of the curve~$\gamma$. Then
there exists a bounded, connected, strongly
pseudoconvex domain~$\Omega$ with
$C^\infty$~regular boundary such that $\Omega$~is 
contained in $G_1\cup
G_2\cup  V$, and $\Omega$ coincides
with $G_1\cup G_2$ outside~$V$.
\end{barbell}

When $G_1$ and~$G_2$ are balls of equal size, and
$\gamma$~is the shortest line segment joining
them, then the domain~$\Omega$ is a ``barbell,''
or dumbbell.

The $C^\infty$~regularity can be changed
everywhere in the statement of the lemma
to $C^k$~regularity, where $k$~is any integer
greater than or equal to~$2$.

\begin{proof}[Proof of the theorem]
I have not claimed that the bounded pseudoconvex
open sets which fail to be Lu Qi-Keng form an open
set in either of the metrics $\rho_1$ or~$\rho_2$,
and I do not know whether or not this is the case
for sets with irregular boundaries. However, if
$\Omega$~is, for example, a 
bounded strongly pseudoconvex open set
with $C^\infty$~regular boundary, and the Bergman kernel
function of~$\Omega$ has zeroes on some connected
component~$\Omega_0$, 
then there is a
$\rho_1$~neighborhood of~$\Omega$ containing no 
pseudoconvex Lu Qi-Keng open set. Indeed, if a
sequence of pseudoconvex open sets converges
to~$\Omega$ in the metric~$\rho_1$, then the
corresponding Bergman kernel functions converge
on~$\Omega_0$ to
the Bergman kernel function of~$\Omega_0$ by the
stability lemma, and by Hurwitz's theorem these
approximating Bergman kernel functions cannot all
be zero-free on~$\Omega_0$. The analogous statement holds for
the metric~$\rho_2$ if one restricts attention to
sets of uniformly bounded Euclidean diameter.

Accordingly, to prove the theorem it will suffice
to construct, arbitrarily close (according to
either $\rho_1$ or~$\rho_2$) to a given bounded
pseudoconvex open set~$G$  a bounded
strongly pseudoconvex open set~$\Omega$ with
$C^\infty$~regular boundary whose Bergman kernel
function does have zeroes on some connected
component; if $G$~is connected, then $\Omega$~should
be connected too.

It is standard that the pseudoconvex open 
set~$G$ can be exhausted from
inside by strongly pseudoconvex open sets with
$C^\infty$~regular boundaries: namely, by sublevel
sets of a smooth, strictly plurisubharmonic exhaustion
function. 
It is evident that
these interior approximating sets converge to~$G$
in both of the metrics $\rho_1$ and~$\rho_2$.
Consequently, there is no loss of generality in
supposing from the start that $G$~is a bounded 
strongly pseudoconvex open set with
$C^\infty$~regular boundary.

Place close to~$G$ a strongly pseudoconvex
domain~$D$ with $C^\infty$~regular boundary and
small diameter, the Bergman kernel function of~$D$
having zeroes. (In~$\C^1$, the domain~$D$ could be
an annulus; in higher dimensions, $D$~could be a
small homothetic copy of the counterexample domain
that I constructed in~\cite{Boas}.)  Then $G\cup
D$ will be a disconnected strongly pseudoconvex
open set that is close to~$G$ in both of the
metrics $\rho_1$ and~$\rho_2$. This open set
$G\cup D$ will serve as the required~$\Omega$ to
prove parts (1) and~(3) of the theorem.

To prove parts (2) and~(4) of the theorem, I need
to produce a connected~$\Omega$ when $G$~is
connected.  To do this, join $G$ to~$D$ with a
closed line segment~$L$, and use the barbell lemma
to construct a sequence of bounded, connected,
strongly pseudoconvex open sets~$\Omega_k$ with
$C^\infty$~regular boundaries, the~$\Omega_k$
being contained in $G\cup D\cup V_k$, where the
$V_k$ are shrinking neighborhoods of the line
segment~$L$. The $\Omega_k$ converge to $G\cup D$
in the metric~$\rho_2$, so the stability lemma and
Hurwitz's theorem imply that the Bergman kernel
function of~$\Omega_k$ has zeroes (on~$D$) when
$k$~is sufficiently large. Since 
the Euclidean distance of 
$D\cup V_k$ from~$G$ is small, $\Omega_k$~is close
to~$G$ in the metric~$\rho_1$ as well as in the
metric~$\rho_2$.  Thus one of the~$\Omega_k$
serves as the required~$\Omega$.
\end{proof}

\begin{proof}[Proof of the stability lemma]
The main point is to prove an $L^2$~approximation theorem
for holomorphic functions. I claim that if $f$~is
a square-integrable holomorphic function on~$U$,
and if a positive~$\epsilon$ is prescribed, then
for all sufficiently large~$j$ there exists a
square-integrable holomorphic function~$g_j$
on~$\Omega_j$ such that
$\|f-g_j\|_{L^2(\Omega_j\cap U)}<\epsilon$ and
$\|g_j\|_{L^2(\Omega_j \setminus U)}<\epsilon$.

I first need to show that the holomorphic functions in
the Sobolev space $W^1(U)$ of square-integrable functions with
square-integrable first derivatives are dense in
the space of square-integrable holomorphic
functions on~$U$. 
This is a consequence of Kohn's global regularity
theorem \cite{Kohn} for the $\dbar$-Neumann
problem with weights.
Namely, for a suitably large positive
number~$t$, the weighted $\dbar$-Neumann
operator~$N_t$ for~$U$ is a bounded operator on
the Sobolev space $W^2(U)$. Consequently, the
corresponding weighted Bergman projection operator~$P_t$,
which satisfies the relation
$P_t=\operatorname{Id}-\dbar_t^* N_t \dbar$, maps
$W^3(U)$ into  the holomorphic subspace of
$W^1(U)$. Now if $f$~is a square-integrable
holomorphic function in~$U$, take a sequence
$\{v_j\}$ of $C^\infty$ functions converging
to~$f$ in $L^2(U)$, and project these functions
by~$P_t$. The functions $P_tv_j$ are holomorphic
functions in $W^1(U)$ that converge to~$f$ in
$L^2(U)$. 

Therefore, there is no loss of
generality in assuming from the start that the holomorphic
function~$f$ lies in
$W^1(U)$. Consequently, $f$~is the restriction
to~$U$ of a function~$F\in W^1(\C^n)$.

It follows from the hypothesis of the lemma that
there is an open neighborhood~$V$ of the closure
of~$U$ such that the $2n$-dimensional Lebesgue 
measure of $V\cap(\Omega_j\setminus U)$ tends to
zero as $j\to\infty$. There is no harm in cutting
off the function~$F$ so that its support lies
inside~$V$. 

The one-form $\dbar F$ is then defined on all
of~$\C^n$, zero on~$U$, zero outside~$V$, and
square-integrable. Since the measure of 
$V\cap(\Omega_j\setminus U)$ shrinks to zero, the
$L^2(\Omega_j)$~norm of~$\dbar F$  tends to zero as
$j\to\infty$. Use H\"ormander's $L^2$~theory
\cite{Hormander} to
solve the equation $\dbar u_j=\dbar F$
on~$\Omega_j$ for a square-integrable
function~$u_j$ whose $L^2(\Omega_j)$ norm is bounded by a constant
(depending only on the uniform bound on the
diameters of the~$\Omega_j$) times the
$L^2(\Omega_j)$ norm
of~$\dbar F$. Thus the norm of~$u_j$ on~$\Omega_j$
tends to zero as $j\to\infty$. Consequently, the
function $g_j:=F-u_j$, which is holomorphic and
square-integrable
on~$\Omega_j$, has norm on~$\Omega_j\cap U$ close to the norm
of~$f$ when $j$~is large. Also, the norm of~$g_j$
on $\Omega_j\setminus U$ tends to zero with the
measure of $V\cap(\Omega_j\setminus U)$. This confirms the
claimed approximation property.

The remainder of the proof of the stability lemma follows
standard lines. However, I mention that I am
dispensing with the hypothesis of monotonicity of
the domains that is typically assumed
\cite[pp.~180--182]{Jarnicki Pflug},
\cite[pp.~36--39]{Skwarczynski}.

Fix a point~$z$ in~$U$. The Bergman kernel
function $K(\cdot,z)$ (for~$U$, or equivalently
for~$\Omega$ when the free variable is in~$U$) is
the unique square-integrable holomorphic
function~$f$ on~$U$ that maximizes $f(z)$ subject
to the nonlinear constraint $f(z)\ge \|f\|^2_{L^2(U)}$. Let
$f_j$~denote the corresponding extremal function
for the approximating domain~$\Omega_j$. By the
mean-value property of holomorphic functions,
$f_j(z)$ is bounded by a constant times $\|f_j\|_{L^2(\Omega_j)}$
times an inverse power of the distance from~$z$ to
the boundary of~$\Omega_j$; the extremal property
of~$f_j$ then implies that $\|f_j\|_{L^2(\Omega_j)}$ too is
bounded by a constant times an inverse power of
the distance from~$z$ to the boundary
of~$\Omega_j$. Therefore the $\|f_j\|_{L^2(\Omega_j)}$ are
uniformly bounded, and so the $f_j$~form a normal
family on~$U$.  Consequently, the $f_j$~have a
subsequence that converges uniformly on compact
subsets of~$U$ to a holomorphic limit
function~$f_\infty$. (Once I show that the
limit~$f_\infty$ actually is~$f$, it will follow
that the original sequence~$\{f_j\}$, not just a
subsequence, converges to~$f$.)

By Fatou's lemma, it follows that the limit
function~$f_\infty$ satisfies
$f_\infty(z)\ge\|f_\infty\|^2_{L^2(U)}$. By the
approximation property proved above, there exists
a square-integrable holomorphic function~$g_j$
on~$\Omega_j$ such that $g_j(z)\ge
\|g_j\|^2_{L^2(\Omega_j)}$, and $g_j(z)\ge
(1-\delta_j)f(z)$, where the positive
numbers~$\delta_j$ tend to zero as $j\to\infty$.
The extremal function~$f_j$ therefore has the
property that $f_j(z)\ge(1-\delta_j)f(z)$.
Consequently, $f_\infty(z)\ge f(z)$. The
uniqueness of the extremal function implies that
$f_\infty=f$. This proves that the Bergman kernel
functions $K_j(w,z)$ for the~$\Omega_j$ converge
pointwise to $K(w,z)$ on $U\times U$.

Since $|K_j(w,z)|^2\le K_j(w,w)K_j(z,z)$ by the
Cauchy-Schwarz inequality, and the right-hand side
is bounded by a constant depending only on the
distances of $z$ and~$w$ from the boundary
of~$\Omega_j$, the functions $K_j(\cdot\,,\cdot)$
form a normal family in $U\times U$. From the
normality and the pointwise convergence just
proved, it is immediate that the convergence is
uniform on compact subsets of $U\times U$.
\end{proof}

\begin{proof}[Proof of the barbell lemma]
In the complex plane~$\C^1$, there is nothing to
prove, for every planar domain is strongly
pseudoconvex. In higher dimensions,
there is no loss of generality in supposing that
the curve~$\gamma$ meets the boundaries of $G_1$
and~$G_2$ transversely, since the barbell~$\Omega$
is not prescribed inside the neighborhood~$V$. 
By \cite[Theorem~4]{Fornaess Zame} (a result that
the authors attribute to \cite{Fornaess Stout}), the set
$\overline G_1\cup \overline G_2 \cup\gamma$ has a
basis of Stein neighborhoods, so there exists a
connected, strongly pseudoconvex domain with
$C^\infty$~regular boundary that outside~$V$ is a
small perturbation of $G_1\cup G_2$.
This conclusion is already enough for the application to the
proof of the main theorem.

The stronger statement that one can find a barbell
that actually matches $G_1\cup G_2$ outside a
neighborhood of the curve~$\gamma$ was
demonstrated by Shcherbina for the case when 
$G_1$ and~$G_2$ are balls \cite[Lemma~1.2
and its Corollary]{Shcherbina}. The
general case follows from this special one because
any strongly pseudoconvex domain can be perturbed
in an arbitrarily small neighborhood of a boundary
point to obtain a new strongly pseudoconvex domain
whose boundary near that point is a piece of the
boundary of a
ball. This can be seen from the patching lemma for strictly
plurisubharmonic functions in
\cite[Lemma~3.2.2]{Eliashberg}
by taking the totally real set there to be a
single point.
\end{proof}

\section*{Open questions}
\begin{enumerate}
\item In the stability lemma, the $C^\infty$
regularity hypothesis can be reduced to
$C^2$~regularity by inspecting Kohn's proof
\cite{Kohn} to see that $C^{k+1}$ boundary
regularity suffices for $W^k$~regularity of the weighted
$\dbar$-Neumann operator; one also needs
techniques as in \cite{Boas Straube} to see that
the weighted Bergman projection has the same
regularity as the weighted $\dbar$-Neumann
operator. Can the hypothesis in the stability
lemma be reduced to
$C^1$~boundary regularity?
\item The conclusion of the theorem---that most
pseudoconvex domains are not Lu Qi-Keng
domains---changes if the topology on domains is
changed. For example, 
any small $C^\infty$~perturbation of the unit ball is a Lu
Qi-Keng domain
\cite{Greene Krantz}.
Does the set of bounded pseudoconvex Lu Qi-Keng
domains with $C^1$~regular boundary have nonempty
interior in the $C^1$~topology on pseudoconvex domains? This is
the case for domains in the complex plane~$\C^1$.

\item My proof of the stability lemma for the
Bergman kernel function uses pseudoconvexity.
Can the word ``pseudoconvex'' be removed from
the statement of the main theorem? 

\item Is every bounded \emph{convex} domain  a
Lu Qi-Keng domain?
\end{enumerate}

\end{document}